\documentclass[12pt]{article}
\usepackage{latexsym, amssymb, amsmath, amscd, amsfonts, epsfig, graphicx, colordvi,verbatim,ifpdf,extarrows}
\usepackage{amsfonts, amsmath, amssymb}
\usepackage{amssymb,amsfonts,amsmath,latexsym,epsfig,cite, psfrag,eepic,color}
\usepackage{amscd,graphics}
\usepackage{latexsym, amssymb,  amsmath,amscd, amsfonts, epsfig, graphicx, colordvi,amsthm}

\usepackage{graphicx}
\usepackage{epstopdf}
\usepackage{color}
\usepackage{ifpdf}
\usepackage{fancybox}
\usepackage[font=small,labelfont=bf,labelsep=none]{caption}
\usepackage{float}

\allowdisplaybreaks

\newtheorem{thm}{Theorem}[section]

\newtheorem{conj}[thm]{Conjecture}

\newtheorem{defi}[thm]{Definition}
\newtheorem{lem}[thm]{Lemma}
\newtheorem{core}[thm]{Corollary}

\def\pf{\noindent{\it Proof.} }
\def\qed{\nopagebreak\hfill{\rule{4pt}{7pt}}
\medbreak}

\numberwithin{equation}{section}

\def\qed{\nopagebreak\hfill{\rule{4pt}{7pt}}
\medbreak}

\newlength{\boxedparwidth}
  {\begin{center} \begin{tabular}{|@{\hspace{.315in}}c@{\hspace{.15in}}|}
                  \hline \\ \begin{minipage}[t]{\boxedparwidth}
                  \setlength{\parindent}{.25in}}%
  {\end{minipage} \\ \\ \hline \end{tabular} \end{center}}

\begin{document}

\begin{center}
{\Large \bf Combinatorial interpretations of truncated versions of a identity of Gauss}
\end{center}

\begin{center}
 {Thomas Y. He}$^{1}$ and
  {S.Y. Liu}$^{2}$ \vskip 2mm

$^{1,2}$ School of Mathematical Sciences, Sichuan Normal University, Chengdu 610066, P.R. China

   \vskip 2mm

  $^1$heyao@sicnu.edu.cn,  $^2$liushuyu@stu.sicnu.edu.cn
\end{center}

\vskip 6mm   {\noindent \bf Abstract.} In 2012, Andrews and Merca proved a truncated theorem on Euler's pentagonal number theorem, which opened up a new study on truncated theta series. In particular, some truncated versions of a identity of Gauss have been proved. In this article, we provide new combinatorial interpretations of the truncated versions of the identity of Gauss in terms of the minimal excludant non-overlined part of an overpartition.



\section{Introduction}

A partition $\pi$ of a positive integer $n$ is a finite non-increasing sequence of positive integers $\pi=(\pi_1,\pi_2,\ldots,\pi_m)$ such that $\pi_1+\pi_2+\cdots+\pi_m=n$. The empty sequence forms the only partition of zero. The $\pi_i$ are called the parts of $\pi$. Let $p(n)$ denote the number of partitions of $n$. The generating function of $p(n)$ was given by Euler.
\[\sum_{n=0}^\infty p(n)q^n=\frac{1}{(q;q)_\infty}.\]
Here and in the sequel, we assume that $|q|<1$ and employ the standard notation:
\[(a;q)_\infty=\prod_{i=0}^\infty(1-aq^i),\]
\[(a;q)_n=\frac{(a;q)_\infty}{(aq^n;q)_\infty},\]
and
\[{M\brack N}=\left\{\begin{array}{ll}\frac{(q;q)_M}{(q;q)_N(q;q)_{M-N}},&\text{if }M\geq N\geq 0,\\
0,&\text{otherwise.}
\end{array}
\right.\]

One of the well-known theorems in the partition theory is Euler's pentagonal number theorem:
\begin{equation}\label{Euler-pentagonal}
(q;q)_\infty=\sum_{n=0}^\infty(-1)^nq^{n(3n+1)/2}(1-q^{2n+1}).
\end{equation}

In 2012, Andrews and Merca \cite{Andrews-Merca-2012} considered a truncated version of \eqref{Euler-pentagonal} and obtained that for $k\geq 1$,
\[\frac{1}{(q;q)_\infty}\sum_{j=0}^{k-1}(-1)^jq^{j(3j+1)/2}(1-q^{2j+1})=1+(-1)^{k-1}\sum_{n=1}^\infty\frac{q^{k(k-1)/2+(k+1)n}}{(q;q)_n}{{n-1}\brack{k-1}},\]
from which they deduced the following theorem for the partition function $p(n)$.

\begin{thm}\cite[Theorem 1.1]{Andrews-Merca-2012}
For $n,k\geq 1$,
\[(-1)^{k-1}\sum_{j=0}^{k-1}(-1)^j\left(p(n-j(3j+1)/2)-p(n-j(3j+5)/2-1)\right)=M_{k}(n),\]
where $M_{k}(n)$ is the number of partitions of $n$ in which $k$ is the least integer that is not a part and there are more parts $>k$ than there are $<k$.
\end{thm}

Apart from Euler's pentagonal number theorem, there is one other classical theta identity due to Gauss \cite[(2.2.12)]{Andrews-1976}:
\begin{equation}\label{Gauss}
1+2\sum_{j=1}^\infty(-1)^jq^{j^2}=\frac{(q;q)_\infty}{(-q;q)_\infty}.
\end{equation}

Inspired by Andrews and Merca's work, Guo and Zeng \cite{Guo-Zeng-2013} established the truncated version of \eqref{Gauss} in 2013. They proved that for $k\geq 1$,
\begin{equation}\label{Guo-Zeng-truncated-Gauss}
\frac{(-q;q)_\infty}{(q;q)_\infty}\left(1+2\sum_{j=1}^k(-1)^jq^{j^2}\right)=1+(-1)^k\sum_{n=k+1}^\infty\frac{(-q;q)_k(-1;q)_{n-k}q^{(k+1)n}}{(q;q)_n}{{n-1}\brack{k}}.
\end{equation}

An overpartition, introduced by Corteel and Lovejoy \cite{Corteel-Lovejoy-2004},  is a partition such that the last occurrence of a number may be overlined.  For example, there are  fourteen overpartitions of $4$.
\[(4),(\overline{4}),(3,1),(\overline{3},1),(3,\overline{1}),(\overline{3},\overline{1}),(2,2),(2,\overline{2}),
\]
\[(2,1,1),(\overline{2},1,1),(2,1,\overline{1}),(\overline{2},1,\overline{1}),(1,1,1,1),(1,1,1,\overline{1}).\]

We impose the following order on the parts of an overpartition.
\begin{equation}\label{order}
\overline{1}<1<\overline{2}<2<\cdots.
\end{equation}

Let $\overline{p}(n)$ be the number of overpartitions of $n$. The generating function of $\overline{p}(n)$ is
\[\sum_{n=0}^\infty\overline{p}(n)q^n=\frac{(-q;q)_\infty}{(q;q)_\infty}.\]

The following corollary is an immediate consequence of \eqref{Guo-Zeng-truncated-Gauss}.

\begin{core}\cite[Corollary 1.2]{Guo-Zeng-2013}
For $n,k\geq 1$,
\begin{equation}\label{Guo-Zeng-truncated-Gauss-1}
(-1)^{k}\left(\overline{p}(n)+2\sum_{j=1}^k(-1)^j\overline{p}(n-j^2)\right)\geq 0,
\end{equation}
with strict inequality if $n\geq (k+1)^2$. For example,
\begin{equation}\label{k=1-case}
\overline{p}(n)-2\overline{p}(n-1)\leq 0,\text{ with strict inequality if }n\geq 4.
\end{equation}
\end{core}

A combinatorial proof of \eqref{k=1-case} was given by Guo and Zeng \cite{Guo-Zeng-2013}. In 2018, Andrews and Merca \cite{Andrews-Merca-2018} provided the following truncated version of \eqref{Gauss}. For $k\geq 1$,
\begin{equation*}\label{proof-use-1}
\frac{(-q;q)_\infty}{(q;q)_\infty}\left(1+2\sum_{j=1}^k(-1)^jq^{j^2}\right)=1+2(-1)^k\frac{(-q;q)_k}{(q;q)_k}\sum_{j=0}^\infty\frac{q^{(k+1)(k+j+1)}(-q^{k+j+2};q)_\infty}{(q^{k+j+1};q)_\infty},
\end{equation*}
from which they deduced the following combinatorial interpretation of  the sum in \eqref{Guo-Zeng-truncated-Gauss-1}. The combinatorial proof of the following theorem was given by Ballantine, Merca, Passary and Yee \cite{Ballantine-Merca-Passary-Yee-2018}.
\begin{thm}\cite[Corollary 8]{Andrews-Merca-2018}\label{Andrews-Merca-truncated-gauss}
For $n,k\geq 1$,
\begin{equation*}
(-1)^{k}\left(\overline{p}(n)+2\sum_{j=1}^k(-1)^j\overline{p}(n-j^2)\right)=\overline{M}_k(n),
\end{equation*}
where $\overline{M}_k(n)$ is the number of overpartitions of $n$ in which the first part larger than $k$ appears at least $k+1$ times.
\end{thm}

The minimal excludant of a partition was introduced by Grabner and Knopfmacher \cite{Grabner-Knopfmacher-2006} under the name ``smallest gap". Recently, Andrews and Newman \cite{Andrews-Newman-2019} undertook a combinatorial study of the minimal excludant of a partition. The minimal excludant of a partition $\pi$ is the smallest positive integer that is not a part of $\pi$, denoted  $mex(\pi)$.
In \cite{Ballantine-Merca-2021}, Ballantine and Merca proved some inequalities involving the minimal
excludant of a partition. They also posed the following conjecture.
\begin{conj}\cite[Conjecture 1]{Ballantine-Merca-2021}\label{conjecturen-Ballantine-Merc}
For $n,k\geq 1$,
\begin{equation*}
\sum_{j=-\infty}^\infty(-1)^j\overline{M}_k(n-j(3j-1))\geq 0,
\end{equation*}
with strict inequality if $n\geq (k+1)^2$.
\end{conj}

Conjecture \ref{conjecturen-Ballantine-Merc} was settled by Yao in \cite{Yao-2025}. Furthermore, Yao showed that for $\ell\geq 1$,
\begin{equation}\label{eqn-Yao-2025}
\sum_{n=0}^\infty\sum_{j=-\infty}^\infty(-1)^j\overline{M}_k(n-\ell j(3j-1)/2)q^n=2\frac{(q^\ell;q^\ell)_\infty}{(q;q)_\infty}\sum_{j=0}^\infty \frac{q^{(k+2j+1)^2}(1-q^{2k+4j+3})}{(q;q^2)_\infty}.
\end{equation}

In \cite{Guo-Zeng-2013}, Guo and Zeng conjectured a stronger inequality than \eqref{Guo-Zeng-truncated-Gauss-1}.

\begin{conj}\cite[(6.4)]{Guo-Zeng-2013}\label{conjecturen-(k+1)}
For $n,k\geq 1$,
\begin{equation}\label{Guo-Zeng-truncated-Gauss-conjuecture}
(-1)^{k-1}\left(\overline{p}(n)+2\sum_{j=1}^k(-1)^j\overline{p}(n-j^2)\right)+\overline{p}(n-k^2)\geq 0,
\end{equation}
with strict inequality if $n\geq k^2$.
\end{conj}

Conjecture \ref{conjecturen-(k+1)} was proved independently by Mao \cite{Mao-2015} and Yee \cite{Yee-2015} and then reconfirmed later by Wang and Yee \cite{Wang-Yee-2019}. In \cite{Xia-Yee-Zhao-2022}, Xia, Yee and Zhao proved the following inequality which implies \eqref{Guo-Zeng-truncated-Gauss-conjuecture}.
\begin{equation*}
(-1)^{k-1}\left(\overline{p}(n)+2\sum_{j=1}^k(-1)^j\overline{p}(n-j^2)\right)+\overline{p}(n-k(k+1))\geq 0.
\end{equation*}

In \cite{Li-2023}, Li obtained the following truncated version of \eqref{Gauss}. For $k\geq 1$,
\begin{equation}\label{proof-use-2}
\frac{(-q;q)_\infty}{(q;q)_\infty}\sum_{j=-k}^{k-1}(-1)^jq^{j^2}=1+(-1)^{k-1}\frac{(-q;q)_k}{(q;q)_{k-1}}\sum_{j=0}^\infty\frac{q^{k(k+j)}(-q^{k+j+1};q)_\infty}{(q^{k+j+1};q)_\infty},
\end{equation}
based on which Li \cite{Li-2023} presented a combinatorial interpretation of the sum in \eqref{Guo-Zeng-truncated-Gauss-conjuecture}.
\begin{thm}\cite[Theorem 1.2]{Li-2023}\label{combinatorial-N-k}
For $n,k\geq 1$,
\[(-1)^{k-1}\sum_{j=-k}^{k-1}(-1)^j\overline{p}(n-j^2)=\overline{N}_k(n),\]
where $\overline{N}_k(n)$ is the number of overpartitions of $n$ in which a part of size $k$ has to be overlined, the smallest part $>k-1$ appears exactly $k$ times and it cannot be overlined.
\end{thm}

In \cite{Wang-Xiao-2024,Wang-Xiao-2025}, Wang and Xiao also obtained \eqref{proof-use-2} and gave a combinatorial interpretation of the sum in \eqref{Guo-Zeng-truncated-Gauss-conjuecture}. Li \cite{Li-2023} generalized the sums in \eqref{Guo-Zeng-truncated-Gauss-1} and \eqref{Guo-Zeng-truncated-Gauss-conjuecture}.

 \begin{thm}\cite[Corollary 5.3]{Li-2023}
 For integers $m$ and $k$ with $m\leq k$,
 \begin{equation}\label{main-eqn-m-k}
 (-1)^{\text{min}\{|m|,k\}}\sum_{j=m}^{k}(-1)^j\overline{p}(n-j^2)\geq 0.
 \end{equation}
 \end{thm}

\section{Main results}

In this section, we will list the results obtained in this article. One of the main objectives of this article is to present a new combinatorial proof of  \eqref{k=1-case}.

There is a generalization of the minimal excludant of a partition. In \cite{Andrews-Newman-2020}, Andrews and Newman defined $mex_{A,a}(\pi)$ to be the smallest positive integer congruent to $a$ modulo $A$ that is not a part of $\pi$.
In \cite{Yang-Zhou-2024}, Yang and Zhou extended the definition of $mex_{A,a}(\pi)$ to overpartitions.

\begin{defi}\cite[Definition 4.1]{Yang-Zhou-2024}
Let $\overline{mex}_{A,a}(\pi)$ be the smallest positive integer congruent to $a$ modulo $A$ that is not a non-overlined part in the overpartition $\pi$. For $n\geq 1$, define
$op_{A,a}(n)$ (resp. $\overline{op}_{A,a}(n)$) to be the number of overpartitions $\pi$ of $n$ such that $\overline{mex}_{A,a}(\pi)$ is congruent to $a$ (resp. $a+A$) modulo $2A$.
\end{defi}

By definition, it is easy to get that for $n\geq 1$,
\begin{equation}\label{a-A-P}
op_{A,a}(n)+\overline{op}_{A,a}(n)=\overline{p}(n).
\end{equation}

For example, we have
\begin{center}
\begin{tabular}{|c|c|c|c|c|c|c|c|c|}
  \hline
  $\pi$&$(4)$&$(\overline{4})$&$(3,1)$&$(\overline{3},1)$&$(3,\overline{1})$&$(\overline{3},\overline{1})$&$(2,2)$&$(2,\overline{2})$\\
  \hline
  $\overline{mex}_{2,1}(\pi)$&$1$&$1$&$5$&$3$&$1$&$1$&$1$&$1$\\
  \hline
\end{tabular}
\end{center}

\begin{center}
\begin{tabular}{|c|c|c|c|c|c|c|}
  \hline
   $\pi$&$(2,1,1)$&$(\overline{2},1,1)$&$(2,1,\overline{1})$&$(\overline{2},1,\overline{1})$&$(1,1,1,1)$&$(1,1,1,\overline{1})$\\
\hline
  $\overline{mex}_{2,1}(\pi)$&$3$&$3$&$3$&$3$&$3$&$3$\\
  \hline
\end{tabular}
\end{center}
Hence, $op_{2,1}(4)=7$ and $\overline{op}_{2,1}(4)=7$.

Yang and Zhou \cite{Yang-Zhou-2024} obtained the following theorem.
\begin{thm}\cite[Theorem 4.2]{Yang-Zhou-2024}\label{useful-theorem-003}
For $n\geq 1$, we have
\begin{equation}\label{useful-theorem-001-1}
op_{2,1}(n)=\overline{p}(n)/{2}.
\end{equation}
\end{thm}

For $n\geq 1$ and $k\geq 0$, define $op_{2,1}(n,k)$ to be the number of overpartitions $\pi$ of $n$ such that $\overline{mex}_{2,1}(\pi)\geq 2k+1$ and $\overline{mex}_{2,1}(\pi)\equiv 2k+1\pmod{4}$. By definition, we have
\begin{equation*}\label{relation-2-1-j}
op_{2,1}(n)=op_{2,1}(n,0)\text{ and }\overline{op}_{2,1}(n)=op_{2,1}(n,1).
\end{equation*}
Combining with \eqref{a-A-P} and \eqref{useful-theorem-001-1}, we can get the following theorem.
\begin{thm}\label{useful-lemma-2}
For $n\geq 1$,
\begin{equation*}
op_{2,1}(n,1)=\overline{p}(n)/{2}.
\end{equation*}
\end{thm}

In this article, we will give another proof of Theorem \ref{useful-lemma-2}. Then, we will present a combinatorial interpretation of the sum in \eqref{main-eqn-m-k} from the point of view of the minimal excludant non-overlined part of an overpartition.
\begin{thm}\label{combiatorial-proof-main-thm-m-k}
 For integers $m$ and $k$ with $m\leq k$,
  \begin{equation}\label{combiatorial-proof-main-eqn-m-k-1}
 (-1)^{\text{min}\{|m|,k\}}\sum_{j=m}^{k}(-1)^j\overline{p}(n-j^2)=op_{2,1}(n,a)+(-1)^{m+k}op_{2,1}(n,b+1)\text{ if }mk> 0,
 \end{equation}
 and
  \begin{equation}\label{combiatorial-proof-main-eqn-m-k-2}
 (-1)^{\text{min}\{|m|,k\}}\sum_{j=m}^{k}(-1)^j\overline{p}(n-j^2)=op_{2,1}(n,a+1)+(-1)^{m+k}op_{2,1}(n,b+1)\text{ if }mk\leq 0,
 \end{equation}
 where   $a=\text{min}\{|m|,|k|\}$ and $b=\text{max}\{|m|,|k|\}$.
\end{thm}

Setting $m=-k<0$ and $m=-k+1\leq 0$ in \eqref{combiatorial-proof-main-eqn-m-k-2}, we can get new combinatorial interpretations of  the sums in \eqref{Guo-Zeng-truncated-Gauss-1} and \eqref{Guo-Zeng-truncated-Gauss-conjuecture}  respectively.

\begin{core}\label{Andrews-Merca-truncated-gauss-new}
For $n,k\geq 1$,
\begin{equation}\label{main-eqn-1}
(-1)^{k}\left(\overline{p}(n)+2\sum_{j=1}^k(-1)^j\overline{p}(n-j^2)\right)=2op_{2,1}(n,k+1),
\end{equation}
and
\begin{equation*}\label{main-eqn-2}
(-1)^{k-1}\left(\overline{p}(n)+2\sum_{j=1}^k(-1)^j\overline{p}(n-j^2)\right)+\overline{p}(n-k^2)=op_{2,1}(n,k)-op_{2,1}(n,k+1).
\end{equation*}
\end{core}

For $k\geq 1$, the generating function of $op_{2,1}(n,k+1)$ is
\begin{equation}\label{gen-op}
\begin{split}
\sum_{n=1}^{\infty}op_{2,1}(n,k+1)&=(-q;q)_\infty\sum_{j=0}^\infty \frac{q^{1+3+\cdots+(2k+4j+1)}(1-q^{2k+4j+3})}{(q;q)_\infty}\\
&=\frac{(-q;q)_\infty}{(q;q)_\infty}\sum_{j=0}^\infty q^{(k+2j+1)^2}(1-q^{2k+4j+3}).
\end{split}
\end{equation}
Then, the following truncated version of \eqref{Gauss} immediately follows from \eqref{main-eqn-1} and \eqref{gen-op}.
\begin{core}For $k\geq 1$,
\[
\frac{(-q;q)_\infty}{(q;q)_\infty}\left(1+2\sum_{j=1}^k(-1)^jq^{j^2}\right)=1+2(-1)^k\frac{(-q;q)_\infty}{(q;q)_\infty}\sum_{j=0}^\infty q^{(k+2j+1)^2}(1-q^{2k+4j+3}).
\]
\end{core}

Combining Theorem \ref{Andrews-Merca-truncated-gauss}, Theorem \ref{combinatorial-N-k} and Corollary \ref{Andrews-Merca-truncated-gauss-new}, we can get the following corollary.
\begin{core}\label{add-coro}For $n,k\geq 1$,
\begin{equation}\label{core-eqn-1}
op_{2,1}(n,k+1)=\overline{M}_k(n)/{2},
\end{equation}
and
\begin{equation*}\label{core-eqn-2}
op_{2,1}(n,k)-op_{2,1}(n,k+1)=\overline{N}_k(n).
\end{equation*}
\end{core}

Then, we can get \eqref{eqn-Yao-2025} by using \eqref{Euler-pentagonal} with $q\rightarrow q^\ell$, \eqref{gen-op}, \eqref{core-eqn-1}, and Euler's partition identity \cite[(1.2.5)]{Andrews-1976}:
\begin{equation*}
\frac{1}{(q;q^2)_\infty}=(-q;q)_\infty.
\end{equation*}

Note that $\overline{M}_0(n)=\overline{p}(n)$ for $n\geq 1$, then by Theorem \ref{useful-lemma-2}, we have
  \[op_{2,1}(n,1)=\overline{M}_0(n)/{2}.\]
  Combining with Corollary \ref{add-coro}, we get the following result.
\begin{core}\label{third-result}For $k\geq 1$,
\begin{equation*}\label{third-eqn}
\sum_{n=1}^\infty\left(\overline{M}_{k-1}(n)-\overline{M}_{k}(n)\right)q^n=2\frac{(-q;q)_k}{(q;q)_{k-1}}\sum_{j=0}^\infty\frac{q^{k(k+j)}(-q^{k+j+1};q)_\infty}{(q^{k+j+1};q)_\infty}.
\end{equation*}
\end{core}

This article is organized as follows. We first give a combinatorial proof of \eqref{k=1-case} in Section 3. Then, we show Theorems \ref{useful-lemma-2} and \ref{combiatorial-proof-main-thm-m-k} in Section 4. Finally, we give an analytic proof of Corollary \ref{third-result} in Section 5.

\section{Combinatorial proof of \eqref{k=1-case}}

Clearly, in order to prove \eqref{k=1-case}, it is equivalent to showing that for $n\geq 1$,
\begin{equation}\label{k=1-case-equiv}
\overline{p}(n)/2\leq \overline{p}(n-1),\text{ with strict inequality if }n\geq 4.
\end{equation}
To do this, we introduce three sets of overpartitions. Bear in mind that the parts in an overpartition are ordered as in \eqref{order}.

For $n\geq 1$, let $\mathcal{A}(n)$ be the set of overpartitions of $n$ such that the smallest part is non-overlined. For example, we have
\[\mathcal{A}(4)=\{(4),(3,1),(\overline{3},1),(2,2),(2,1,1),(\overline{2},1,1),(1,1,1,1)\}.\]

For an overpartition $\pi$ in $\mathcal{A}(n)$, if we change the smallest part to an overlined part, then we get an overpartition of $n$ such that the smallest part is overlined, and vice versa. This implies that the number of overpartitions in $\mathcal{A}(n)$ is $\overline{p}(n)/2$.

For $n\geq 1$, let  $\mathcal{B}(n-1)$ be the set of overpartitions $\lambda$ of $n-1$ such that if $\overline{1}$ appears in $\lambda$ then the smallest part lager than $1$ in $\lambda$  is greater than or equal to two plus the number of parts $1$ in $\lambda$. For example, we have
\[\mathcal{B}(3)=\{(3),(\overline{3}),(2,1),(\overline{2},1),(2,\overline{1}),(1,1,1),(1,1,\overline{1})\}.\]

For $n\geq 1$, we set $\mathcal{C}(n-1)$ be the set of overpartitions of $n-1$  not in $\mathcal{B}(n-1)$, that is, $\mathcal{C}(n-1)$ is the set of overpartitions $\lambda$ of $n-1$ such that $\overline{1}$ appears in $\lambda$ and the smallest part larger than $1$ in $\lambda$ is less than two plus the number of parts $1$ in $\lambda$. For example, we have
\[\mathcal{C}(3)=\{(\overline{2},\overline{1})\}.\]

To give a combinatorial proof of \eqref{k=1-case-equiv}, it suffices to show the following theorem.
\begin{thm}
There exists a bijection between $\mathcal{A}(n)$ and $\mathcal{B}(n-1)$ for $n\geq 1$ and $\mathcal{C}(n-1)$ is nonempty for $n\geq 4$.
\end{thm}

\pf  For $n\geq 1$, let $\pi=(\pi_1,\ldots,\pi_{m-1},\pi_m)$ be an overpartition in  $\mathcal{A}(n)$. By definition, we know that $\pi_m$ is a non-overlined part, which implies that $\overline{1}$ does not appear in $\pi$. Assume that  $\pi_m=t$, we consider the following two cases.

Case 1: $t=1$. In this case, we set $\lambda=(\pi_1,\ldots,\pi_{m-1})$.

Case 2: $t\geq 2$. In this case, we set
\[\lambda=(\pi_1,\ldots,\pi_{m-1},\underbrace{1,\ldots,1}_{(t-2)'\text{s}},\overline{1}).\]

In either case, we get an overpartition $\lambda$ in $\mathcal{B}(n-1)$. Obviously, the process above is reversible. Now, we have built a bijection between $\mathcal{A}(n)$ and $\mathcal{B}(n-1)$ for $n\geq 1$.

It remains to show that $\mathcal{C}(n-1)$ is nonempty for $n\geq 4$. It can be checked that for $n\geq 4$,
\[(\overline{2},\underbrace{1,\ldots,1}_{(n-4)'\text{s}},\overline{1})\]
is an overpartition in $\mathcal{C}(n-1)$. This completes the proof. \qed

\section{Proofs of Theorems \ref{useful-lemma-2} and \ref{combiatorial-proof-main-thm-m-k}}

The objective of this section is to give the proofs of Theorems \ref{useful-lemma-2} and \ref{combiatorial-proof-main-thm-m-k}. To do this, we need the following lemma.
\begin{lem}\label{useful-lemma-1}For $n,j\geq 1$,
\begin{equation*}\label{useful-lemma-001}
\overline{p}(n-j^2)=op_{2,1}(n,j)+op_{2,1}(n,j+1).
\end{equation*}
\end{lem}

\pf For an overpartition $\pi$  of $n-j^2$, if we add the parts $1,3,\ldots,2j-1$ into $\pi$, then we get an overpartition $\lambda$ of $n$ with $\overline{mex}_{2,1}(\lambda)\geq 2j+1$, and vice versa. It implies that  $\overline{p}(n-j^2)$ equals the number of overpartitions $\lambda$ of $n$ with $\overline{mex}_{2,1}(\lambda)\geq 2j+1$.

On the other hand, it follows from definition that $op_{2,1}(n,j)+op_{2,1}(n,j+1)$ is also the number of overpartitions  $\lambda$ of $n$ with $\overline{mex}_{2,1}(\lambda)\geq 2j+1$. This completes the proof.  \qed

Then, we give a proof of Theorem \ref{useful-lemma-2}  with the aid of Lemma \ref{useful-lemma-1}.

{\noindent \bf Proof of Theorem \ref{useful-lemma-2}.} Divide both sides of \eqref{Gauss} by $\frac{(q;q)_\infty}{(-q;q)_\infty}$, we can get
\begin{equation}\label{proof-lemma}
\frac{(-q;q)_\infty}{(q;q)_\infty}\left(1+2\sum_{j=1}^\infty(-1)^jq^{j^2}\right)=1.
\end{equation}

Comparing the coefficients of $q^n$ on the both sides of \eqref{proof-lemma}, we have
\[\overline{p}(n)+2\sum_{j=1}^\infty(-1)^j\overline{p}(n-j^2)=0.\]
Combining with Lemma \ref{useful-lemma-1}, we get
\[\overline{p}(n)=-2\sum_{j=1}^\infty(-1)^j\left(op_{2,1}(n,j)+op_{2,1}(n,j+1)\right)=2op_{2,1}(n,1).\]
The proof is complete.  \qed

 We conclude this section with a proof of Theorem \ref{combiatorial-proof-main-thm-m-k}.

{\noindent\bf Proof of Theorem \ref{combiatorial-proof-main-thm-m-k}.}  We first show \eqref{combiatorial-proof-main-eqn-m-k-1}. There are the following two cases.

Case 1: $k\geq m\geq 1$. In this case, we have $\text{min}\{|m|,k\}=|m|=m$, $\text{min}\{|m|,|k|\}=|m|=m$ and $\text{max}\{|m|,|k|\}=|k|=k$. It follows from Lemma \ref{useful-lemma-001} that
\begin{align*}
&\quad(-1)^{\text{min}\{|m|,k\}}\sum_{j=m}^{k}(-1)^j\overline{p}(n-j^2)\\
&=(-1)^m\sum_{j=m}^{k}(-1)^j\left(op_{2,1}(n,j)+op_{2,1}(n,j+1)\right)\\
&=(-1)^m\left((-1)^mop_{2,1}(n,m)+(-1)^kop_{2,1}(n,k+1)\right)\\
&=op_{2,1}(n,m)+(-1)^{m+k}op_{2,1}(n,k+1),
\end{align*}
which agrees with \eqref{combiatorial-proof-main-eqn-m-k-1} for $k\geq m\geq 1$.

Case 2: $m\leq k\leq -1$. In this case, we have $\text{min}\{|m|,k\}=k$, $\text{min}\{|m|,|k|\}=|k|=-k$ and $\text{max}\{|m|,|k|\}=|m|=-m$. Again by  Lemma \ref{useful-lemma-001}, we have
\begin{align*}
&\quad(-1)^{\text{min}\{|m|,k\}}\sum_{j=m}^{k}(-1)^j\overline{p}(n-j^2)\\
&=(-1)^k\sum_{j=-k}^{-m}(-1)^j\overline{p}(n-j^2)\\
&=(-1)^k\sum_{j=-k}^{-m}(-1)^j\left(op_{2,1}(n,j)+op_{2,1}(n,j+1)\right)\\
&=(-1)^k\left((-1)^{-k}op_{2,1}(n,-k)+(-1)^{-m}op_{2,1}(n,-m+1)\right)\\
&=op_{2,1}(n,-k)+(-1)^{m+k}op_{2,1}(n,-m+1).
\end{align*}
We arrive at \eqref{combiatorial-proof-main-eqn-m-k-1} for $m\leq k\leq -1$.

Now, we have completed the proof of \eqref{combiatorial-proof-main-eqn-m-k-1}. Then, we proceed to show \eqref{combiatorial-proof-main-eqn-m-k-2}. Under the condition that $m\leq k$ and $mk\leq 0$, we have $m\leq 0\leq k$. Clearly, we have
\begin{equation}\label{add-eqn-0*}
\sum_{j=m}^{k}(-1)^j\overline{p}(n-j^2)=\sum_{j=m}^{-1}(-1)^j\overline{p}(n-j^2)+\overline{p}(n)+\sum_{j=1}^{k}(-1)^j\overline{p}(n-j^2).
\end{equation}

 In view of \eqref{combiatorial-proof-main-eqn-m-k-1}, we get
\begin{equation}\label{add-eqn-*}
\sum_{j=m}^{-1}(-1)^j\overline{p}(n-j^2)=-\left(op_{2,1}(n,1)+(-1)^{m+1}op_{2,1}(n,-m+1)\right),
\end{equation}
and
\begin{equation}\label{add-eqn-**}
\sum_{j=1}^{k}(-1)^j\overline{p}(n-j^2)=-\left(op_{2,1}(n,1)+(-1)^{k+1}op_{2,1}(n,k+1)\right).
\end{equation}

Substituting \eqref{add-eqn-*} and \eqref{add-eqn-**} into \eqref{add-eqn-0*} and combining with Theorem \ref{useful-lemma-2}, we have
\begin{equation}\label{proof-add***}
\sum_{j=m}^{k}(-1)^j\overline{p}(n-j^2)=(-1)^{m}op_{2,1}(n,-m+1)+(-1)^{k}op_{2,1}(n,k+1).
\end{equation}

Then, we consider the following two cases.

Case 1: $k\geq -m\geq 0$. In this case, we have $\text{min}\{|m|,k\}=|m|=-m$, $\text{min}\{|m|,|k|\}=|m|=-m$ and $\text{max}\{|m|,|k|\}=|k|=k$. Appealing to \eqref{proof-add***}, we get
\[\quad(-1)^{\text{min}\{|m|,k\}}\sum_{j=m}^{k}(-1)^j\overline{p}(n-j^2)=op_{2,1}(n,-m+1)+(-1)^{m+k}op_{2,1}(n,k+1).\]
So, \eqref{combiatorial-proof-main-eqn-m-k-2} is valid for $k\geq -m\geq 0$.

Case 2: $-m>k\geq 0$. In this case, we have $\text{min}\{|m|,k\}=k$, $\text{min}\{|m|,|k|\}=|k|=k$ and $\text{max}\{|m|,|k|\}=|m|=-m$. Using \eqref{proof-add***}, we get
\[\quad(-1)^{\text{min}\{|m|,k\}}\sum_{j=m}^{k}(-1)^j\overline{p}(n-j^2)=op_{2,1}(n,k+1)+(-1)^{m+k}op_{2,1}(n,-m+1),\]
and thus \eqref{combiatorial-proof-main-eqn-m-k-2} is satisfied for $-m>k\geq 0$. This completes the proof.  \qed

\section{Analytic proof of Corollary \ref{third-result}}

In this section, we aim to give an analytic proof of  Corollary \ref{third-result}. Note that for $k\geq 0$, the generating function of $\overline{M}_k(n)$ is
\[\sum_{n=1}^{\infty}\overline{M}_k(n)q^n=2\frac{(-q;q)_k}{(q;q)_k}\sum_{j=0}^\infty\frac{q^{(k+1)(k+j+1)}(-q^{k+j+2};q)_\infty}{(q^{k+j+1};q)_\infty},\]
so it remains to show
\begin{equation}\label{finally-M-K}
\begin{split}
&\quad\frac{(-q;q)_{k-1}}{(q;q)_{k-1}}\sum_{j=0}^\infty\frac{q^{k(k+j)}(-q^{k+j+1};q)_\infty}{(q^{k+j};q)_\infty}-\frac{(-q;q)_k}{(q;q)_k}\sum_{j=0}^\infty\frac{q^{(k+1)(k+j+1)}(-q^{k+j+2};q)_\infty}{(q^{k+j+1};q)_\infty}\\
&=\frac{(-q;q)_k}{(q;q)_{k-1}}\sum_{j=0}^\infty\frac{q^{k(k+j)}(-q^{k+j+1};q)_\infty}{(q^{k+j+1};q)_\infty}.
\end{split}
\end{equation}

Dividing both sides of \eqref{finally-M-K} by $\frac{(-q;q)_{k-1}}{(q;q)_{k}}$, we find that it is equivalent to
\begin{equation}\label{finally-M-K-0}
\begin{split}
&\quad(1-q^{k})\sum_{j=0}^\infty\frac{q^{k(k+j)}(-q^{k+j+1};q)_\infty}{(q^{k+j};q)_\infty}-(1+q^k)\sum_{j=0}^\infty\frac{q^{(k+1)(k+j+1)}(-q^{k+j+2};q)_\infty}{(q^{k+j+1};q)_\infty}\\
&=(1-q^{2k})\sum_{j=0}^\infty\frac{q^{k(k+j)}(-q^{k+j+1};q)_\infty}{(q^{k+j+1};q)_\infty}.
\end{split}
\end{equation}

Then, we proceed to show  \eqref{finally-M-K-0}.

{\noindent \bf Proof of \eqref{finally-M-K-0}.} Note that
\begin{align*}
&\quad(1+q^k)q^{(k+1)(k+j+1)}\\
&=q^{(k+1)(k+j+1)}+q^{(k+1)(k+j+1)+k}\\
&=q^{k(k+j+1)}(q^{k+j+1}-1+1)+q^{k(k+j+2)}(q^{k+j+1}+1-1)\\
&=q^{k(k+j+1)}(1-q^k)-q^{k(k+j+1)}(1-q^{k+j+1})+q^{k(k+j+2)}(1+q^{k+j+1}),
\end{align*}
so we get
\begin{align*}
&\quad(1+q^k)\sum_{j=0}^\infty\frac{q^{(k+1)(k+j+1)}(-q^{k+j+2};q)_\infty}{(q^{k+j+1};q)_\infty}\\
&=\sum_{j=0}^\infty\left((1-q^k)\frac{q^{k(k+j+1)}(-q^{k+j+2};q)_\infty}{(q^{k+j+1};q)_\infty}-\frac{q^{k(k+j+1)}(-q^{k+j+2};q)_\infty}{(q^{k+j+2};q)_\infty}\right.\\
&\qquad\qquad\qquad\left.+\frac{q^{k(k+j+2)}(-q^{k+j+1};q)_\infty}{(q^{k+j+1};q)_\infty}\right)\\
&=(1-q^k)\sum_{j=1}^\infty\frac{q^{k(k+j)}(-q^{k+j+1};q)_\infty}{(q^{k+j};q)_\infty}-\sum_{j=1}^\infty\frac{q^{k(k+j)}(-q^{k+j+1};q)_\infty}{(q^{k+j+1};q)_\infty}\\
&\qquad\qquad\qquad+\sum_{j=0}^\infty\frac{q^{k(k+j+2)}(-q^{k+j+1};q)_\infty}{(q^{k+j+1};q)_\infty}.
\end{align*}

Then, we have
\begin{align*}
&\quad(1-q^{k})\sum_{j=0}^\infty\frac{q^{k(k+j)}(-q^{k+j+1};q)_\infty}{(q^{k+j};q)_\infty}-(1+q^k)\sum_{j=0}^\infty\frac{q^{(k+1)(k+j+1)}(-q^{k+j+2};q)_\infty}{(q^{k+j+1};q)_\infty}\\
&=(1-q^{k})\frac{q^{k^2}(-q^{k+1};q)_\infty}{(q^{k};q)_\infty}+\sum_{j=1}^\infty\frac{q^{k(k+j)}(-q^{k+j+1};q)_\infty}{(q^{k+j+1};q)_\infty}-\sum_{j=0}^\infty\frac{q^{k(k+j+2)}(-q^{k+j+1};q)_\infty}{(q^{k+j+1};q)_\infty}\\
&=\sum_{j=0}^\infty\frac{q^{k(k+j)}(-q^{k+j+1};q)_\infty}{(q^{k+j+1};q)_\infty}-\sum_{j=0}^\infty\frac{q^{k(k+j+2)}(-q^{k+j+1};q)_\infty}{(q^{k+j+1};q)_\infty}\\
&=(1-q^{2k})\sum_{j=0}^\infty\frac{q^{k(k+j)}(-q^{k+j+1};q)_\infty}{(q^{k+j+1};q)_\infty}.
\end{align*}
We arrive at \eqref{finally-M-K-0}. This completes the proof.  \qed

\end{document}